\documentclass[12pt,a4paper]{article}
\usepackage{amsmath,amssymb,amsthm,amsfonts,amscd,euscript,verbatim}
\usepackage{t1enc}
\usepackage[cp1251]{inputenc}

\theoremstyle{plain}
\newtheorem{theo}{Theorem}[section]

\newtheorem{pr}[theo]{Proposition}

\theoremstyle{remark}

\theoremstyle{remark}
\newtheorem{defi}[theo]{Definition}
\newtheorem*{notat}{Notation}

 \newcommand\de{\Delta}

\newcommand\ok{{\mathfrak{O}_K}}

\newcommand\ol{{\mathfrak{O}_L}}

\newcommand\ovl{{{\overline L}}}

\newcommand\gs{\mathfrak{FGS}{}}

\newcommand\gm{\mathfrak{M}}

\newcommand\ob{^{-1}} \newcommand\lan{\langle}
\newcommand\ra{\rangle}
\newcommand\nde{\triangleleft}

\newcommand\ns{\{0\}}
\newcommand\bff{\mathbf{f}}
\newcommand\bv{\mathbf{V}}

\DeclareMathOperator\homm{\operatorname{Hom}}
\DeclareMathOperator\ext{\operatorname{Ext}}

\DeclareMathOperator\ke{\operatorname{Ker}}
\DeclareMathOperator\cok{\operatorname{Coker}}

\DeclareMathOperator\car{\operatorname{Cart}}

\DeclareMathOperator\cl{\operatorname{Cl}}
\DeclareMathOperator\spe{\operatorname{Spec}}

\begin{document}
 \title{Finite flat commutative group schemes over complete discrete valuation
 rings:  classification in terms of Cartier modules, structural results; application to reduction of Abelian varieties}
 \author{M.V. Bondarko}
 \date{ 13.11.2004}
\begin{abstract} This paper is a summary of author's results on
finite flat commutative group schemes. The properties of the
generic fibre functor are discussed. A complete classification  of
finite local flat commutative group schemes over mixed
characteristic complete discrete valuation rings in terms of their
Cartier modules (defined by Oort) is given. We also state several
properties of tangent space of these schemes. These results are
applied to the study of reduction of Abelian varieties. A finite
$p$-adic semistable reduction criterion is formulated. It looks
especially nice for the ordinary reduction case. The plans of the
proofs are described.

Keywords: finite group scheme, Cartier module, tangent space,
formal group, Abelian variety, semistable reduction,  local field.

MSC 2000: {14L15, 14L05, 14G20, 11G10, 11S31}.

\end{abstract}
 \maketitle

\begin{notat}
$K$ is a mixed characteristic complete discrete valuation  field
with residue field of characteristic $p$,
 $L$ is a finite  extension of
$K$; $\mathfrak{O}_K\subset\mathfrak{O}_L$ are their rings of
integers, $e$ is the absolute ramification index of $L$,
$s=[\log_p(pe/(p-1))]$, $e_0=[L:(K_{nr}\cap L)]$ ($e_0=e(L/K)$ in the perfect residue field
case), $l'=s+v_p(e_0)+1$, $l=2s+v_p(e_0)+1$; $\ovl$ denotes the
residue field of $\ol$;  $\gm$ is the maximal ideal of $\ol$;
$\pi\in \gm$ is some uniformizing element of $L$.

In this paper a 'group scheme' will (by default) mean a finite
flat commutative group scheme, $S/\ol$ means a finite flat
commutative group scheme over $\ol$.

For finite group schemes $S,T$ we  write $S\nde T$  if $S$ is a
closed subgroup scheme of $T$.
\end{notat}

\section{The category of finite flat commutative group schemes; the generic fibre results}

We denote by $\gs_R$ the category of finite flat commutative
$p$-group schemes (i.e. annihilated by a power of $p$) over a base
ring $R$.

The goal of this paper is the study of $\gs_\ol$ and of Abelian
varieties over $L$. We note that a certain classification of
$\gs_\ol$ for $\ovl$ being perfect was given by Breuil (see
\cite{clb}); yet  that classification is inconvenient for several
types of problems.

It is well known that any finite flat group scheme over $L$ is
\'etale; hence $\gs_L$ is equivalent to the category of finite
modules over the absolute Galois group of $L$. In particular, this
category is abelian.

Hence it is natural to consider the generic fibre functor $GF:S\to
S_L=S\times_{\spe\ol}\spe L$. $GF$ is  faithful and defines a
one-to-one correspondence between closed subgroup schemes of $S$
and closed subgroup schemes of $S_L$ (see \cite{re}).

It was also proved by Raynaud  in the case $e<p-1$ that $GF$ is
full; besides $\gs_\ol$ is an Abelian category. Neither of this
facts is true for larger values of $e$. Moreover, one cannot apply
Raynaud's methods in the case $e>p-1$.

Yet the following important result is valid.

\begin{theo}\label{maingf}
If $S,T/\ol$ are group schemes, $g:S_L\to T_L$ is an $L$-group scheme morphism,
 then there exists
an $h:S\to T$ over $\ol$ such that $h_L=p^s g$.

\end{theo}

Note that $s=0$ for $e<p-1$; therefore Theorem \ref{maingf}
generalizes the fullness result of Raynaud.

Hence $GF$ is 'almost full'. One easily checks that the result is
sharp, i.e. the value of $s$ is the best possible.

Theorem \ref{maingf} also can be considered as a finite analogue
of fullness of the generic fibre functor for $p$-divisible groups
(proved by Tate). Besides, it implies Tate's result (see
\cite{ta}) immediately.

The main tool of the proof is the Cartier module functor for
finite local group schemes. It will be defined below.

A similar statement for $\ext^1$ follows easily from Theorem \ref{maingf} and the
Cartier module theory for group schemes.

\section{Formal groups; Cartier modules}

Our basic method is resolving finite group schemes by means of
$p$-divisible groups (in particular, by finite height formal group
laws).

We remind the  Cartier module theory for formal group laws. Here
we describe a modified version that was used in \cite{0} and
\cite{01} (cf. \cite{1} and \cite{z}).

We denote by $C$ the category  of additive subgroups of $L[[\de]]^m,\ m>0$.
For $C_1,C_2\in C$, $\dim C_i=m_i$ we define
 $$C(C_1,C_2)=A\in M_{m_2\times m_1}\mathfrak{O}_K:\ AC_1\subset C_2.$$

 For
 $h\in (L[[\Delta]])^m$, the coefficients of $h $ are equal to
$h_{i}=\sum_{l\ge 0} c_{il}\Delta^l$, $c_{il}\in L$, we
define
\begin{equation}\label{ap}
h(x)=(h_i(x)),\ 1\le i\le m,\text{ where }h_i(x)=\sum_{l\ge 0}c_{il}
x^{p^l}.
\end{equation}

For  an $m$-dimensional formal group law $F/\ol$ we consider $D_F=\{f\in
L[[\de]]^m: \exp_F(f(x))\in \ol[[x]]^m\}$, where $\exp_F\in L[[X]]^m$
 is the composition inverse to the logarithm of $F$.
 In particular, for
$m=1$ we have $\sum a_i\de^i\in D_F\iff\exp_F(\sum a_ix^{p^i})\in
\ol[[x]]$.

Then the Cartier theory easily implies the following fact.

\begin{pr}\label{mcart1}
1. $F\to D_F$ defines a full  embedding of the  category of formal
groups over $\ol$ into $C$.

 If $f:F_1\to F_2$, $f\equiv AX\mod \deg 2$, $A\in M_{m_2\times m_1}\ol$,
 then the associated map  $f_*:D_{F_1}\to D_{F_2}$ is the  multiplication by $A$.

2. $D_{F_1}=D_{F_2}$ if and only if the groups $F_1$ and $F_2$ are
strictly isomorphic, i.e. there exists an isomorphism whose linear
term is given by the identity matrix.\end{pr}

Now we briefly remind the notion of the Cartier ring.

For a commutative ring $Q$ and a $Q$-algebra $P$ one can introduce
the following operators on $P[[\de]]$.

For $f=\sum_{i\ge 0} c_i\de^i \in P[[\de]],\ a\in Q$  we define
 $$\bv f=f \de;\ \bff f=\sum_{i>0} pc_i \de^{i-1};\ \lan a\ra f=\sum a^{p^i}c_i\de^i .$$

 $\car(Q)$ (the Cartier $p$-ring, see \cite{1}) is the ring that is generated by $\bv,\bff,\lan a\ra,\ a\in Q$ satisfying
 certain natural relations (see \cite{1}, section 16.2, \cite{02}, and \cite{04}).

If
 $M$ is a $\car(Q)$-module, then $M/\bv M$
has a natural structure of a $Q$-module defined via $a\cdot
(x\mod\bv M) =\lan a\ra x \mod\bv M $ for any  $x\in M$.

 We introduce an important definition (see \cite{02}).

\begin{defi}
1. For $\car$-modules $M\subset N$ we write $M\nde N$, if for any
$x\in M, \bv x\in N$, we have $x\in  M$. We call $M$ a {\bf closed}
submodule of  $N$.

2. $\car$-module $N$ is called {\bf separated} if $\ns\nde N$,
i.e. $N$ has no $\bv$-torsion.
\end{defi}

We denote $\car(\ol)$ by $\car$; we call $\car$-modules Cartier
modules.

Note that we don't define closed subsets of Cartier modules. Yet
we could define a topology on any $\car$-module $M$ whose closed
subsets would be $C\subset M:\ \bv x\in C\implies x\in C$. Then
any  $\car$-module homomorphism would be a continuous map.

We define the closure of a subset $N$ of a Cartier module $M$
as the smallest closed $\car$-submodule of $M$ that contains $N$.

\begin{pr}\label{mcart2}
1. For any $F/\ol$ the group $D_F$ is a $\car$-module via the
action of operators defined above; it is canonically
$\car$-isomorphic to the module of $p$-typical curves for $F$ (see
\cite{1} and \cite{z}).

2. $C(D_{F_1}, D_{F_2})=\car(D_{F_1}, D_{F_2})$.

3. $M\in C$ is equal to $D_F$ for some $m$-dimensional $F/\ol$ iff $C\nde
L[[\de]]^m$ and $C\mod\de=\ol^m$.
\end{pr}

In the papers \cite{01} and \cite{0} two functors on the category
of formal groups were defined. The first (called the fraction
part) was similar to certain defined by Grothendieck, Messing and
Fontaine;
 yet it was defined
in a quite different way and was described more precisely than the
functor in the book  \cite{2}. The behaviour of the fraction part
is (in some sense) linear.

The second functor (denoted by $M_F$) described the obstacle for the Fontaine's
functor to be an embedding of categories. For a finite height
formal group $F$  the value of $M_F$ can be described by means of $D_{F_\pi}$. Here
$F_\pi=\pi\ob(\pi X,\pi Y)$. Since the coefficients of $F_\pi$
tend to $0$ quickly, the obstacle functor is 'finite'. One may say
that its complexity is killed by $\bff^s$ (see Proposition 3.2.2
of \cite{02}).

\section{Cartier-Oort modules of local group schemes}

Let $S$ be a local group scheme over $\ol$; let $0\to S\to F\to
G\to 0$ be its resolution by means of finite height formal groups.
We define $C(S)=\cok(D_F\to D_G)$.

In the paper \cite{oo} it was proved that
 $S\to C(S)$ is a well defined functor on the category of local (finite flat commutative) group schemes over $\ol$; it defines an embedding of this category into the category of $\car$-modules.

We call $C(S)$ the Oort module of $S$. The theory of Oort also can be used
when the base ring is a field of characteristic $p$. In this case $\bff$ corresponds
to the Frobenius, $\bv$ corresponds to the Verschiebung operator
(see \cite{z}).

Now we state the main classification result. It completely
describes the properties of the Oort functor.

\begin{theo}\label{clsc}

I1. Closed submodules of $C(S)$ are in one-to-one correspondence
with closed subgroup schemes of $S$.

2. If $M=C(S)$, $N=C(H)\nde M$, where $H\nde S$, then $M/N\approx
C(S/H)$.

3. Conversely, exact sequences (as fppf-sheaves, i.e. the inclusion
is a closed embedding) of local schemes induce exact sequences of
Oort modules.

II If $f:S\to T$ is a local group scheme morphism, then $\ke f_*=C(\ke
f)$, where $f_*$ is the induced Oort modules homomorphism; we
consider the kernel in the category of flat group schemes.

III  A $\car$-module $M$ is isomorphic to $C(S)$ for $S$ being a
finite flat commutative local group scheme over $\ol$ if and only
if $M$ satisfies the following conditions.

1. $M/\bv M$ is a finite length $\ol$-module.

2. $M$ is separated.

3. $\cap_{i\ge 0} \bv^i M=\{0\}$.

4. $M=\cl_M(\lan\pi\ra M)$.

IV The minimal dimension of a finite height formal group $F$ such
that $S$ can be embedded into $F$ is equal to $\dim_\ol(C(S)/\bv
C(S))$ (i.e. to  the number of indecomposable $\ol$-summands of
$C(S)/\bv C(S)$).

V $M=C(\ke [p^r]_F)$ for an $m$-dimensional formal group $F$ if
and only if in addition to the conditions of III, $p^r M=0$ and
$M/\bv M\approx (\ol/p^r\ol)^m$.

VI If $S,T$ are local, then $\ext^1(S,T)=\ext^1_{\car}(C(S),
C(T))$. Here we consider extensions in the category of finite flat
group schemes, whence the definition of an exact sequence is the
same as in part I3.

\end{theo}

We introduce a natural definition of the tangent space $TS$ for
a finite group scheme $S$.

\begin{defi}\label{tgs}
For a finite flat groups scheme $S/\ol$ we denote by $TS$ the
$\ol$-dual of $J/J^2$ (i.e. $\homm_\ol (J/J^2,L/\ol)$), where $J$
is the augmentation ideal of the affine algebra of $S$.
\end{defi}

It is well known that the tangent space of a group scheme is equal
(i.e. naturally isomorphic) to the tangent space of its local
part. Besides, if $P$ is any (unitial commutative) $\ol$-algebra
then the (suitably defined) tangent space of $S_P=S\times_{\spe
\ol} \spe P$ is canonically isomorphic to $TS\otimes_ \ol P$.

We state the main properties of the tangent space functor.

\begin{theo}\label{dimgs}

I $TS$ is naturally isomorphic to $C(S_0)/\bv C(S_0)$, where $S_0$
is the local part of $S$.

II $f:S\to T$ is a closed embedding of local group schemes if and
only if the induced map on the tangent spaces is an embedding.

III If $0\to H\to S\to T\to 0$ is an exact sequence of local group
schemes (in the category of fppf-sheaves, i.e. $H\nde S$) then the
corresponding sequence of tangent spaces is also exact.

IV For a local group scheme $S$ the following numbers are equal.

1. The $\ol$-dimension of $J/J^2$.

2. The $\ol$-dimension of $C(S)/\bv C(S)$.

3. The minimal dimension of a finite height formal group $F$ such
that $S\nde F$.

V A local group scheme $S$ is equal to $\ke [p^r]_F$ for some
$m$-dimensional finite height formal group $F/\ol$ if and only if
$p^r S=0$, and $TS\approx (\ol/p^r\ol)^m$.
\end{theo}

\section{Finite  criteria for reduction of Abelian varieties}

As an application of the results on finite group schemes
 certain finite $p$-adic criteria for semistable
and ordinary reduction of Abelian varieties were proved. We call
these criteria  finite because in contrast to Grothendieck's
criteria (see \cite{gro}) it is sufficient to check certain
conditions on some finite $p$-torsion subgroups of $V$ (instead
of the whole $p$-torsion).

We recall that an Abelian variety (over $\ok$ or $\ol$) is called
an ordinary reduction one (or just ordinary) if   the connected
component of $0$ of the reduction of the N\'eron model of $A$ is
an extension of a torus by an ordinary Abelian variety (over
$\ovl$).

In particular, an ordinary variety has semistable reduction. It
can be easily seen that an Abelian variety is ordinary  if the
formal group of its N\'eron module is finite height and of
multiplicative type.

For example, a semistable reduction elliptic curve is either
ordinary or supersingular.

Let $V$ be  an Abelian variety  of dimension $m$ over $K$
 that has semistable reduction over $L$.

\begin{theo}\label{mared}

I $V$ has semistable reduction over $K$ if and only if
  for there exists a finite flat  group scheme $H/\ok$ such that
 $TH_\ol\supset (\ol/p^{l}\ol)^m$ (i.e. there exists an embedding)
and there exists a monomorphism
 $g:H_K\to\ke[p^{l}]_{V,K}$.

II $V$ has ordinary  reduction over $K$ if and only if
for some  $H_K\subset \ke[p^{l}]_{V,K}$ and $M$ unramified
  over $K$ we have
  $H_M\cong (\mu_{p^{l},M})^m$.
Here $\mu$ denotes the group scheme of roots of unity.
 \end{theo}

Part I is a vast generalization of Theorem 5.3 of \cite{co} where
the case $e<p-1$, $V$ of good reduction over $\ol$, was
considered.

Finite $l$-adic criteria (see \cite{sz1}) seem to be easier to
use; yet they don't allow to check whether the reduction is
ordinary.

If the reduction of $A$ over $L$ is good then $l$ can be replaced
by $l'$.

\section{Short plans for the proofs of the main statements}

Proposition \ref{mcart2} is an easy implication of the usual
Cartier theory. It easily implies parts I and II of Theorem
\ref{clsc}, and parts I -- III of Theorem \ref{dimgs}. The proof
the necessity of conditions of part III in Theorem \ref{clsc} is
also more-or-less easy.

To prove sufficiency of conditions of part III in Theorem \ref{clsc}
one applies the explicit description of the Cartier module of a
formal group (see section 27.7 of \cite{1}) and constructs a
formal group $F$ such that $M$ is a $\car$-factor of $D_F$.
A formal group of dimension $\dim_\ol(M/\bv M)$ can be chosen. Next
one proves that a finite height formal group can be chosen. In
this case $M$ will be equal to $D_F/N$ for some $N\nde D_F$ such that $N$ is
$\car$-isomorphic to $D_G$ for a finite height formal group $G$.
Lastly one verifies that $M=C(S)$ for $S$ being the kernel of a
certain isogeny $h:G\to F$. Under the conditions of part V of
Theorem \ref{clsc} (and Theorem \ref{dimgs}) one obtains that
$G\approx F$, and $h=\ke[p^r]_F$.

Parts IV and V of Theorem \ref{dimgs} are reformulations of the corresponding parts of Theorem \ref{clsc} in terms of tangent spaces.
Part VI of Theorem \ref{clsc}
follows from the fact the the conditions of part III are preserved
by extensions.

Now we sketch the proof of Theorem \ref{maingf}.
First the following important results on the
reductions of group schemes are proved.

\begin{pr}\label{red}

1. If the map $h:S\to T$ of $\ol$-group schemes is injective on
the generic fibre, then the kernel of the reduction map (as a
kernel of a group scheme morphism over $\ovl$) is annihilated by
$\bff^s$.

2. If the map $h:S\to T$ of $\ol$-group schemes is surjective on
the generic fibre, then the cokernel of the reduction map  is
annihilated by $\bv^s$.

In the imperfect residue field case we extend $\ovl$ so that
$Fr^{-s}(\cok \overline{h})$ will be defined over $\ovl$.
\end{pr}

Part 1 is proved by an analysis of the properties of the
'obstacle' functor (see the end of \S2); part 2 follows immediately
as the dual of part 1.

In the proof of Proposition \ref{red} (and in several other
places) the local-\'etale exact sequence for finite flat
commutative group schemes is used to reduce the problem to the
study of local group schemes (and hence to formal groups). Some of
these reduction reasonings are quite complicated.

Next one proceeds to the proof of Theorem \ref{maingf} using the
fact $\bff^s\bv^s=p^s$.

For the proofs of reduction criteria explicit Cartier module descent is used.

In \cite{02} the following important statements were  were proved.
For the first an explicit descent reasoning for $D_F$ was used.
The second was proved using flat descent; the details of the proof
are rather complicated.

\begin{pr}\label{oldde}
1. Let $F$ be a finite height formal group over $\ol$. Suppose
that  its generic fibre
 (as a  $p$-divisible group) $F_L=F\times_{\spe \ol} \spe L$ is defined over $K$,
i.e. there exists a $p$-divisible group $Z_K$ over $K$ such that
\begin{equation}\label{cong}
Z_K\times_{\spe K} \spe L \cong F_L.
\end{equation}
Suppose that for $t=v_p(e_0)+1$ and some group scheme $T/\ok$ we have $\ke [p^t]_Z\approx T\times_{\spe \ok} \spe K$.
Suppose also that  this isomorphism combined with the isomorphism  (\ref{cong})
is the generic fibre of a certain isomorphism
 $T \times_{\spe \ok} \spe \ol \cong \ke [p^t]_{F}$.
Then $Z_K\approx F'_K$ for some formal group $F'/\ok$.

2. Let $V$ be a $p$-divisible group over $\ol$. Suppose that  its
generic fibre
  is defined over $K$ and  its local part is  defined over $\ok$
  i.e. there exist $p$-divisible groups  $U/K$ and $S_0/\ok$, an isomorphsim
    $f:U_L\to V_L$, an isomorphism
$h:S_{0\ol}\to V$ of $S_{0\ol}$ with the local part of $V$, and an
imbedding $i:S_{0K}\to U$ such that  $h_L=f_L\circ i_L$.

  Then $V$ is defined over $\ok$ if and only if the inertia group of $K$ acts trivially on
  $U/i(S_{0K})(F)$ where $F$ is the algebraic closure of $K$.

\end{pr}

Using this, Theorem \ref{maingf}, and Cartier modules of group
schemes one can prove a certain good reduction criterion for
Abelian varieties (see \cite{02}). We don't formulate that
criterion here.

Using part 1 of Proposition \ref{oldde}, Theorem \ref{maingf}, and
a certain tangent space argument one proves the following fact.

\begin{pr}\label{degr}
Let $V$ be a $p$-divisible group  over $K$, let $Y$ be a
$p$-divisible group  of dimension $m$ over $\ol$ (i.e. its local
part is a formal group of dimension $m$).
 Suppose that
$V\times_{\spe K} \spe L\cong Y\times_{\spe \ol} \spe L$. We
denote by $G$ the local part of $Y$, by $J$ the corresponding
subgroup of $V$.

Then the following conditions are equivalent:

I There exists a $p$-divisible group $Z$ over $\ok$ such that
$J\cong Z\times_{\spe \ok} \spe K$.

II For some (finite flat commutative) group scheme $H/\ok$ we have
 $TH_\ol\approx (\ol/p^{l'}\ol)^m$ and there exists a monomorphism
 $g:H_K\to\ke[p^{l'}]_{V,K}$.

III We have $TH_\ol\supset (\ol/p^{l'}\ol)^m$ (i.e. there exists
an embedding); there exists a monomorphism
 $g:H_K\to\ke[p^{l'}]_{V,K}$.

\end{pr}

Note that II $\implies$ III is obvious; I $\implies$ II follows
immediately from  part V of Theorem \ref{dimgs}.

The following result on Abelian varieties is useful. By a certain
duality argument the statement is reduced to Grothendieck's
criterion on semistable reduction of Abelian varieties (see
Proposition 5.13 part c in \cite{gro}).

Let $V/K$ denote an Abelian variety that has semistable reduction
over $L$,  denote by $V_f/L$  the formal part of the $p$-torsion
of $V$ (i.e. the part corresponding to the formal group $F$ of the
N\'eron model of $V$ over $\ol$). It is easily seen that  $V_f$ is
equal to $V_{fK}\times_{\spe K}\spe L$ for a certain canonically
defined $p$-divisible group $V_{fK}$ over $K$.

\begin{theo}\label{dual}

$V$ has semistable reduction over $K$ if and only if there exists
a formal $p$-divisible group $G/\ok$ such that the generic fibre
of $G$ (as a $p$-divisible group) is isomorphic to $V_{fK}$.
\end{theo}

Now we sketch the proof of Theorem \ref{mared}.

Let $V_f$  denote the finite part of  $T_p(V)$ (the $p$-torsion of
$V$) considered as a $p$-divisible group over $L$.

If $V$ has semistable reduction over $K$ then
 $V_f$   corresponds to a certain $m$-dimensional
finite height formal group $Y$  defined over $\ok$. Therefore we
can take $H=[\ke p^l]_{Y,\ok}$. It will be of multiplicative type
(i.e. dual-\'etale) if $V$ is ordinary.

For the converse implication a certain tangent space argument
along with Proposition \ref{degr} proves that $V_f$ is defined
over $\ok$. Then Theorem \ref{dual} proves part I.

Lastly it remains to notice that a $p$-divisible group $Y$ over
$\ol$ is of multiplicative type  iff $\ke[p]_Y$ is. Then an easy
tangent space calculation proves part II.


\begin{thebibliography}{1}

 \bibitem{01}
  Bondarko M. V., {Explicit classification of formal groups over complete
discrete valuation fields with imperfect residue field}
(Russian)//
  Trudy St. Peterburgskogo Matematicheskogo Obsh'estva, vol. 11, 2005,
  P. 1--35.



\bibitem{02}
Bondarko M.V., {Finite flat commutative group schemes over
complete discrete valuation fields I: the generic fibre functor, a
 finite wild criterion for good reduction of Abelian varieties}//
  to appear in Mat. Izv. Akad. Nauk, 2006.



 \bibitem{04}
 Bondarko M.V., {Finite flat commutative group schemes over complete discrete valuation
 rings II: classification, tangent spaces, and semistable reduction of Abelian varieties},
 G\"ottingen, 2004.


\bibitem{0} Bondarko M. V., Vostokov S. V.,
{Explicit classification of formal groups over local fields}//
 Proc. Steklov Inst. Math. 2003, no. 2 (241), P. 35--57.

\bibitem{clb}
Breuil C.,  {Groupes $p$-divisibles, groupes finis et modules
 filtr\'es}// Ann. of Math. 2002. vol. 152.  no. 2, P. 489--549.


\bibitem{co}
 Conrad B., {Finite group schemes over bases with low
ramification}// Compositio Mathematica. 1999. v. 119. P. 239--320.


\bibitem{2}%[Fontaine]
 Fontaine J.M., {Groupes $p$-divisibles sur les corps locaux}//
Asterisque. 1977. n. 47--48. Soc. Math. France,  Paris.

\bibitem{gro}
Grothendieck A.,  {S\'eminare de g\'eom\'etrie algebrique 7 I (Expose IX)},
 Lecture Notes in Mathematics, vol. 288.
Springer-Verlag, Berlin--Heidelberg--New York, 1972.
%P. 313--523.


 \bibitem{1}
 Hazewinkel M., {  Formal groups and applications}.
Springer-Verlag, Berlin--Heidelberg--New York, 1978.

 \bibitem{oo}
  Oort F., {Dieudonn\'e modules of
 finite local group schemes}// Indag. Math. v. 37. 1975. P. 103--123.



\bibitem{re}
Raynaud M., {Schemas en groupes de type $(p,\dots,p)$}// Bull. Soc.
Math. France. 1974. v. 102. P. 241--280.


\bibitem{sz1}
Silverberg A., Zarhin Yu. G., {Reduction of abelian varieties}//
The arithmetic and geometry of algebraic cycles (Banff, AB, 1998),
NATO Sci. Ser. C Math. Phys. Sci., 548,  2000. P. 495--513.



\bibitem{ta}
Tate J., {$p$-divisible groups}//  Proc. Conf. Local Fields
(Driebergen, 1966).  P. 158--183 Springer, Berlin, 1967.

\bibitem{z}
Zink Th., {Cartiertheorie kommutativer formaler Gruppen}.
Teubner-Texte zur Mathematik 68,  BSB B. G. Teubner
Verlagsgesellschaft, Leipzig, 1984.

\end{thebibliography}
\end{document}